# Closed Form Evaluations of Some Exponential Sums

N. A. Carella, March 2011.

*Abstract:* This note provides new closed form evaluations of a few classes of exponential sums associated with elliptic curves and hyperelliptic curves.

Keywords: Character sums, Exponential sums, Elliptic curves, Complex multiplication.

Mathematics Subject Classification: Primary 10G05, Secondary 14K22.

## 1 Introduction

This work deals with the closed forms evaluations of several exponential sums. Currently there is a handful classes of exponential sums that can be evaluated in closed forms. Extensive coverage of this subject is given in [BT], see also the research problem 14, p. 497. The exponential sum

$$S(f) = \sum_{x=0}^{p-1} \left( \frac{f(x)}{p} \right) \qquad (1)$$

has been evaluated in closed forms for nine classes of cubic polynomials $f_n(x) \in \mathbf{F}_p[x]$, and a few other polynomials, see [PD], [SK], [JM], et cetera. The nine classes of cubic polynomials are the followings:

(1) $f_1(x) = x^3 + ax$,
(2) $f_2(x) = x(x^2 + 4ax + 2a^2)$,
(3) $f_3(x) = x^3 + a$,
(4) $f_7(x) = x(x^2 + 21ax + 112a^2)$,
(5) $f_{11}(x) = x^3 - 2^5 \cdot 3 \cdot 11a^2x + 2^4 \cdot 7 \cdot 11^2 a^3$,
(6) $f_{19}(x) = x^3 - 2^3 \cdot 19a^2x + 2 \cdot 19^2 a^3$,
(7) $f_{43}(x) = x^3 - 2^4 \cdot 5 \cdot 43a^2x + 2 \cdot 3 \cdot 7 \cdot 43^2 a^3$,
(8) $f_{67}(x) = x^3 - 2^3 \cdot 5 \cdot 11 \cdot 67a^2x + 2 \cdot 7 \cdot 31 \cdot 67^2 a^3$,
(9) $f_{163}(x) = x^3 - 2^4 \cdot 5 \cdot 23 \cdot 29 \cdot 163a^2x + 2 \cdot 7 \cdot 11 \cdot 19 \cdot 127 \cdot 163^2 a^3$,

where $0 \neq a \in \mathbb{Z}$, see Theorem 4. The CM elliptic curve $E: y^2 = f_n(x)$ and its quadratic twisted curve $E': sy^2 = f_n(x)$, where $s$ is a quadratic nonresidue are closely linked to exponential sums $S(f_n)$ in question. Each of these elliptic curve has complex multiplication by the ring of integers $\mathbb{Z}[\sqrt{-n}]$ or $\mathbb{Z}[(1+\sqrt{-n})/2]$ in the quadratic field $\mathbb{Q}(\sqrt{-n})$.

The new evaluations are derived from the nine classes of cubic polynomials and a few other polynomials.

The new nine classes of polynomials derived from $f_n(x)$ are the following.

(1) $g_1(x) = x^4 + a$,
(2) $g_2(x) = x^4 + 4ax^2 + 2a^2$,
(3) $g_3(x) = x^6 + a$,
(4) $g_7(x) = x^4 + 21ax^2 + 112a^2$,
(5) $g_{11}(x) = x^6 - 2^5 \cdot 3 \cdot 11 a^2 x^2 + 2^4 \cdot 7 \cdot 11^2 a^3$,
(6) $g_{19}(x) = x^6 - 2^3 \cdot 19 a^2 x^2 + 2 \cdot 19^2 a^3$,
(7) $g_{43}(x) = x^6 - 2^4 \cdot 5 \cdot 43 a^2 x^2 + 2 \cdot 3 \cdot 7 \cdot 43^2 a^3$,
(8) $g_{67}(x) = x^6 - 2^3 \cdot 5 \cdot 11 \cdot 67 a^2 x^2 + 2 \cdot 7 \cdot 31 \cdot 67^2 a^3$,
(9) $g_{163}(x) = x^6 - 2^4 \cdot 5 \cdot 23 \cdot 29 \cdot 163 a^2 x^2 + 2 \cdot 7 \cdot 11 \cdot 19 \cdot 127 \cdot 163^2 a^3$.

These polynomials are associated with the hyperelliptic curves $E_n : y^2 = g_n(x)$, and the associated exponential sums $S(g_n)$. The new closed form evaluations are the following.

**Theorem 1.** (1) If $n = 1, 2,$ or $7$, then

$$S(g_n) = \sum_{x=0}^{p-1} \left(\frac{g_n(x)}{p}\right) = \begin{cases} A & \text{if } p \text{ is inert,} \\ A + \left(\frac{a}{p}\right)u & \text{if } 4p = u^2 + nv^2, \end{cases} \quad (2)$$

where the integer $u$ is uniquely determined by quadratic symbol relation $(u \mid n) = (2 \mid p)$, and the quadratic exponential sum

$$A = \sum_{x=0}^{p-1}\left(\frac{f_n(x)/x}{p}\right). \quad (3)$$

(2) If $n = 3, 11, 19, 43, 67,$ or $163$, then

$$S(g_n) = \sum_{x=0}^{p-1}\left(\frac{g_n(x)}{p}\right) = \begin{cases} B & \text{if } p \text{ is inert,} \\ B + \left(\frac{a}{p}\right)u & \text{if } 4p = u^2 + nv^2, \end{cases} \quad (4)$$

where the integer $u$ is uniquely determined by $(u \mid n) = (2 \mid p)$, $g_n(x) = (x + a_1)(x + a_2)(x + a_3)(x + a_4) \in \mathbf{F}_p[x]$, and

$$B = -1 - \left(\frac{\alpha}{p}\right)(H(\beta) \bmod p) \quad \text{with} \quad \alpha \equiv (a_1 - a_4)(a_2 - a_3), \; \beta \equiv \frac{(a_1 - a_3)(a_2 - a_4)}{(a_1 - a_4)(a_2 - a_3)} \bmod p, \quad (5)$$

and $H(x)$ is the Hasse invariant.

The proof, which is built up from the simpler closed form evaluations, appears in Section 3.1.

Other widely used families of nonCM elliptic curves of genus $g = 1$ are the following:

(10) $F_\beta(x) = x(x-1)(x-\beta)$, or $F_\beta(x) = (x^2 - 1)(x^2 - \beta)$, where $0, 1 \neq \beta \in \mathbb{Z}$,    Legendre form.
(11) $G_\beta(x) = (k^2 x^2 - 1)(x^2 - \beta)$, $0, \pm 1 \neq \beta \in \mathbb{Z}$, and $k \neq 0$,    Newton form.

(12) $H_d(x) = (x^2 - c^2)(c^2 dx^2 - 1)$, $c, d \in \mathbb{Z}$, and $cd(1 - c^4 d) \neq 0$,                                      Edward form.

The closed form evaluations of these exponential sums are discussed in Section 3.2.

## 2. Preliminary

Let $f(x) \in \mathbf{F}_p[x]$ be a polynomial of degree $deg(f) = d > 0$, and let $\chi(x) = (x \mid p)$ be the quadratic character. Any exponential sum $S(f)$ has a closed form evaluation in terms of the characteristics of the continued fractions $[A_0, A_1, \ldots, A_s]$ and $[a_0, a_1, \ldots, a_t]$ of the rational functions $(f(x)-1)/(x^p - x)$ and $(f(x)+1)/(x^p - x)$ as

$$S(f) = \sum_{x=0}^{p-1} \chi(f(x)) = \sum_{x=0}^{p-1} \left(\frac{f(x)}{p}\right) = \begin{cases} \deg(a_t) & \text{if } n_f = t, \\ \deg(A_s) & \text{if } n_f = s, \end{cases} \qquad (6)$$

see [LN, p. 236] for the exact analysis. However, only a handful of classes of exponential sums $S(f)$ have been evaluated in terms of the characteristics of the coefficients of the polynomials $f(x)$ and the primes $p$.

This section recalls a few basic results to facilitate the proofs. The evaluations of the quadratic character sum, cubic sum, and a few other cases in terms of the characteristics of the polynomials $f(x)$ and the primes $p$.

### 2.1. The Number of Points on A Curve.
The exponential sums $S(f)$ under investigation arise in the calculations of the number of points

$$\#C(\mathbf{F}_p) = \#\{(x, y) : y^2 - f(x) = 0, \text{ and } x, y \in \mathbf{F}_p\} = p + 1 + S(f(x)) \qquad (7)$$

of then algebraic curves $C : y^2 = f(x)$ of genus $g \geq 0$. In terms of the characteristic of the prime $p$, this is written as

$$\#C(\mathbf{F}_p) = 1 + \sum_{x=0}^{p-1}\left(1 + \left(\frac{f(x)}{p}\right)\right) = p + 1 + \sum_{x=0}^{p-1}\left(\frac{f(x)}{p}\right) = p + 1 - \pi - \bar{\pi}, \qquad (8)$$

where the trace $Tr(\pi) = -\pi - \bar{\pi} = S(f)$, and $p = \pi\bar{\pi}$ in the quadratic field $\mathbb{Q}(\sqrt{-n})$.

### 2.2. Simple Exponential Sums and Identities.

**Lemma 2.** Let $t \in \mathbf{F}_p$. Then

$$\sum_{x \in \mathbf{F}_q} x^t \equiv \begin{cases} 0 & \text{if } q-1 \neq \mod t, \\ q-1 & \text{if } q-1 = 0 \mod t. \end{cases} \qquad (9)$$

To confirm this statement, let $\theta$ be a generator of the multiplicative group of $\mathbf{F}_p$, and examine the sum of the geometric progression $1, \theta, \theta^2, \ldots, \theta^{p-2}$. The convention $0^0 \equiv 0 \mod p$ in finite fields of characteristic $p$ is observed.

The linear exponential sum $S(ax + b)$ is identically zero: A simple change of variables, and the fact that the finite field has equal number of quadratic residues and quadratic nonresidues lead to

$$S(f) = \sum_{x=0}^{p-1} \left(\frac{ax+b}{p}\right) = \sum_{x=0}^{p-1} \left(\frac{x}{p}\right) = 0. \tag{10}$$

**Quadratic Exponential Sums**. The quadratic exponential sum $S(ax^2 + bx + c)$ has a closed form evaluation for every polynomial. The analysis is a straight forward algebraic manipulation.

***Lemma 3.*** Let $p > 2$ be a prime and let $f(x) = ax^2 + bx + c \in \mathbf{F}_p[x]$, then

$$S(f) = \sum_{x=0}^{p-1} \left(\frac{ax^2+bx+c}{p}\right) = \left(\frac{a}{p}\right) \begin{cases} -1 & \text{if } b^2 - 4ac \neq 0, \\ p-1 & \text{if } b^2 - 4ac = 0. \end{cases} \tag{11}$$

The proof appears in [LN, 230].

**2.3. Cubic Exponential Sums**. The analysis of the cubic exponential sums dates back about a century. The simplest cases $S(x^3 + a)$ and $S(x^3 + ax)$ were studied by Jacobsthal and Schrutka in the early 1900, see [WS], [LN], and similar references.

***Theorem 4.*** ([PD]) Let $n = 1, 2, 3, 7, 11, 19, 43, 67, 163$, and suppose that $y^2 = f_n(x)$ has good reduction at $p$, exempli gratia, $p \nmid 2an$, then

$$S(f_n) = \sum_{x=0}^{p-1} \left(\frac{f_n(x)}{p}\right) = \begin{cases} 0 & \text{if } p \text{ is inert}, \\ \left(\frac{a}{p}\right)u & \text{if } 4p = u^2 + nv^2, \end{cases} \tag{12}$$

or $p = u^2 + nv^2$, where the integer $u$ is uniquely determined by quadratic symbol relation $(u \mid n) = (2 \mid p)$.

This is the culmination of many decades of works on the cubic exponential sums by various workers. The proofs of various special cases are given several papers, for starter, see [LN], [BT], [PD], [SK], and [MN].

For the simpler cubic exponential sums $S(x^3 + a)$ and $S(x^3 + ax)$, finer analysis are available, see [BT, p 190], [RS].

The elliptic curve $E_n : y^2 = f_n(x)$ has complex multiplication by the full ring of algebraic integers $\mathbb{Z}[\sqrt{-n}]$ or $\mathbb{Z}[(1+\sqrt{-n})/2]$, $n = 1, 2, 3, 7, 11, 19, 43, 67, 163$. The group of points on an algebraic curve $E : y^2 = f(x)$ is defined by $E(\mathbf{F}_p) = \{(x, y) : f(x) - y^2 = 0\}$ and its cardinality is given by $\#E(\mathbf{F}_p) = p + 1 + S(f)$. Similarly, the twisted curve $E' : sy^2 = f(x)$ has the group of points $E'(\mathbf{F}_p) = \{(x, y) : f(x) - sy^2 = 0\}$ and its cardinality is given by $\#E(\mathbf{F}_p) + \#E'(\mathbf{F}_p) = 2p + 2$.

**2.4. Transformation and Reduction Formulae.** Several classes of transformations and degree reduction formulae are known for evaluating character sums of degree $d \geq 3$ in terms of simpler ones of degrees $< d$. Two of these transformation formulae are explored here.

***Lemma 5.*** Let $f(x) \in \mathbf{F}_p[x]$ be a polynomial. Then

$$\sum_{x=0}^{p-1} \left(\frac{f(x^2)}{p}\right) = \sum_{x=0}^{p-1} \left(\frac{xf(x)}{p}\right) + \sum_{x=0}^{p-1} \left(\frac{f(x)}{p}\right). \tag{13}$$

The more general version of this transformation identity is given in [WS]. The degree reduction formula below for evaluating character exponential sums of degree $d > 3$ in terms of simpler ones of degrees $< d$ is given [BT, p. 207]. The special case of degree $d = 4$ is as follows.

**Lemma 6.** Let $p$ be a prime and let $f(x) = (x + a_1)(x + a_2)(x + a_3)(x + a_4)$, then

$$S(f) = \sum_{x=0}^{p-1} \left( \frac{(x+a_1)(x+a_2)(x+a_3)(x+a_4)}{p} \right) = -1 - \left( \frac{\alpha}{p} \right) S(F_\beta), \tag{14}$$

where $F_\beta(x) = x(x - 1)(x - \beta)$, $0, 1 \neq \beta \in \mathbb{Z}$, and the parameters $\alpha$ and $\beta$ are

$$\alpha \equiv (a_1 - a_4)(a_2 - a_3) \bmod p, \quad \text{and} \quad \beta \equiv \frac{(a_1 - a_3)(a_2 - a_4)}{(a_1 - a_4)(a_2 - a_3)} \bmod p. \tag{15}$$

The parameters $\alpha$ and $\beta$ are not independent of the permutation of the roots $a_{\pi(1)}, a_{\pi(2)}, a_{\pi(3)}, a_{\pi(4)}$, where $\pi \in \text{Sym}(n)$ is a permutation of $n$ elements. So the evaluation of the exponential sum $S(f)$ is given up to a permutation of the roots of the polynomial $f(x)$. As the permutation map $\pi \in \text{Sym}(4)$ varies, there are six possible values of each parameter $\alpha \in \{ \pm \alpha_1, \pm \alpha_2, \pm \alpha_3 \}$ and $\beta \in \{ \beta_1^{\pm 1}, \beta_2^{\pm 1}, \beta_3^{\pm 1} \}$, where $0 \neq \alpha_i, \beta_i \in \mathbf{F}_p$.

The polynomial $F_\beta(x) = x(x - 1)(x - \beta)$, $0, 1 \neq \beta \in \mathbb{Z}$, satisfies the identities

(i) $F_\beta(x) = F_{1/\beta}(x)$, (ii) $F_\beta(x) = \left( \frac{-1}{p} \right) F_{1-\beta}(x)$, (iii) $F_{\beta^2}(x) = \left( \frac{\beta}{p} \right) F_{(1+\beta)^2/4\beta}(x)$.

The application of the reduction formula requires information about the reducibility/irreducibility of certain cubic polynomials. The *discriminant test* can be used to help identify the reducibility/irreducibility of any cubic polynomial $f(x) = x^3 + ax^2 + bx + c$ of discriminant $D = a^2b^2 - 4b^3 - 4a^3c - 27c^2 + 18abc$ over the prime field $\mathbf{F}_p$ or an extension of it. The information is extracted from the value of the quadratic symbol test $\left( \frac{D}{p} \right) = (-1)^{s+1}$, where $s$ is the number of irreducible factors in $f(x)$ over $\mathbf{F}_p$. More precise criteria are given in [SP].

**2.5. Hasse Invariant.** Given a prime $p = 2n + 1$, the Hasse polynomial (or Hasse invariant) is defined by

$$H(x) = (-1)^n \sum_{k=0}^{n} \binom{n}{k}^2 x^k, \tag{16}$$

confer [HR, p. 261], [SN, p. 141], etc. The Hasse polynomial is a solution of the differential equation $4x(1-x)\frac{d^2}{dx^2}y + 4(1-2x)\frac{d}{dx}y - 1 = 0$ in the ring of polynomials $\mathbf{F}_q[x]$, see [SN, P. 142]. It is linked to the $n$th Legendre polynomial $P_n(x)$ by the formula $\pm H(x) = (1-x)^n P_n((1+x)/(1-x))$.

**Lemma 7.** The Hasse polynomial $H(x)$ has simple roots in its splitting field $\mathbf{F}_q$ or a quadratic extension $\mathbf{F}_q[\theta]$ of $\mathbf{F}_q$.

**Lemma 8.** ([BM]) For a prime $p > 3$, the following hold:

(i) The number of linear factors of $H(x) \in \mathbf{F}_q[x]$ is $N_1 = \begin{cases} 0 & \text{if } p \equiv 1 \bmod 4, \\ 3h(-p) & \text{if } p \equiv 3 \bmod 4. \end{cases}$

(ii) The number of quadratic factors of $H(x) \in \mathbf{F}_q[x]$ is $N_2 = \begin{cases} h(-p)/2 & \text{if } p \equiv 1 \bmod 4, \\ (3h(-p)-1)/2 & \text{if } p \equiv 3 \bmod 8, \\ (h(-p)-1)/2 & \text{if } p \equiv 7 \bmod 8, \end{cases}$ where $h(-p)$ is the class number of the quadratic field $\mathbb{Q}(\sqrt{-p})$.

The first formula for the number $N_1(p)$ of roots of $H(x) \in \mathbf{F}_p[x]$ in $\mathbf{F}_p$ appeared earlier in [YM]. From these formulas it immediately follows that if the polynomial $H(x) \in \mathbf{F}_p[x]$ does not factor completely in $\mathbf{F}_p$, then it has less than $h(-p) = O(p^{1/2})$ roots in $\mathbf{F}_p$. The Hasse invariant is used to evaluate the character sum $S(F_\beta)$ connected with the polynomial $F_\beta(x) = x(x-1)(x-\beta)$.

**Lemma 9.** Let $F_\beta(x) = x(x-1)(x-\beta)$, $\beta \neq 0, 1 \in \mathbf{F}_p$, $p = 2m+1$. Then

$$S(F_\beta) = \sum_{x=0}^{p-1} \left( \frac{x(x-1)(x-\beta)}{p} \right) \equiv (-1)^m \sum_{k=0}^{m} \binom{m}{k}^2 \beta^k \bmod p \equiv H(\beta) \bmod p. \tag{17}$$

Proof: Compute the expression

$$S(F_\beta) = \sum_{x=0}^{p-1} \left( \frac{x(x-1)(x-\beta)}{p} \right) \equiv \sum_{x=0}^{p-1} (x(x-1)(x-\beta))^{(p-1)/2} \bmod p \equiv H(\beta) \bmod p, \tag{18}$$

refer to [HR, p. 251], and [SN, p. 140 ]. ∎

For $p > 2$, *supersingular* elliptic curves are characterized by the condition $H(\beta) \equiv 0 \bmod p$, see [HR] and [SN] for other equivalent criteria. Elliptic curves with complex multiplication are supersingular for every prime $p \equiv 3 \bmod 4$, and ordinary for every prime $p \equiv 1 \bmod 4$, confer [SN, p. 144].

As an example, Lemma 8 implies that the elliptic curve $E: y^2 = F_\beta(x) = x(x-1)(x-\beta)$, $0, 1 \neq \beta$, is not supersingular for any $\beta \in \mathbf{F}_p$, whenever the prime $p = 4m+1$.

**Theorem 10.** (Manin) Let $E: y^2 = F_\beta(x) = x(x-1)(x-\beta)$, $0, 1 \neq \beta$, be an elliptic curve over the rational numbers $\mathbb{Q}$ with good reduction at a prime $p$, and let $a_p = -S(F_\beta)$ be the trace of Frobenius at $p$, then $H(\beta) \equiv a_p \bmod p$. In particular,

$$\#E(\mathbf{F}_p) = \#\{ (x, y) : y^2 - F_\beta(x) = 0, \text{ and } x, y \in \mathbf{F}_p \} = p + 1 + (H(\beta) \bmod p).$$

More generally, for an elliptic curve $E : y^2 = f(x)$, the number of points satisfies $\#E(\mathbf{F}_q) \equiv 1 - H(\beta)^{(q-1)/(p-1)} \bmod p$.

**2.6. Other Exponential Sums.** Several other related evaluations of character sums are included in this section.

**Definition 11.** If $p - 1 \equiv 0 \bmod k$, and $a \neq 0$, the Jacobson sums are defined by

$$\phi_k(a) = \sum_{x=0}^{p-1} \left(\frac{x}{p}\right)\left(\frac{x^k+a}{p}\right) \quad \text{and} \quad \psi_k(a) = \sum_{x=0}^{p-1}\left(\frac{x^k+a}{p}\right). \tag{19}$$

**Proposition 12.** If $p \equiv k + 1 \bmod 2k$ then $\phi_k(a) = 0$, and $\phi_{2k}(a) = 0$.

Proof : As $(p - 1)/k$ is odd number, the set $\{\, x^k \bmod p : 0 < x < p \,\}$ of residues modulo $p$ is just a permutation of the set $\{\, x^{2k} \bmod p : 0 < x < p \,\}$, see [B, p. 184] for more details. ∎

**Theorem 13.** Let $k > 1$ and let $p$ be a prime such that $p - 1 \equiv 0 \bmod k$, then

(i) $\psi_k(a) = 0$ if $\gcd(k, p-1) = 1$. (ii) $\psi_k(a) = \sum_{x=0}^{p-1}\left(\frac{x^k+a}{p}\right) \equiv -\sum_{i=0}^{k/2}\binom{kf}{2fi}a^{f(k-2i)} \bmod p$ if $\gcd(p-1, k) \neq 1$.

Proof: Take $p = 2kf + 1$. Routine calculations show that

$$\psi_k(a) = \sum_{x=0}^{p-1}\left(\frac{x^k+a}{p}\right) \equiv \sum_{x=0}^{p-1}(x^k+a)^{(p-1)/2} \bmod p \equiv \sum_{x=0}^{p-1}\sum_{d=0}^{kf}\binom{kf}{d}a^{kf-d}x^{kd} \bmod p \tag{20}$$

The inner sum vanishes for $kd \not\equiv 0 \bmod p - 1$, but equals $-1$ for $kd \equiv 0 \bmod p - 1$, see Lemma 2. Hence,

$$\psi_k(a) \equiv -\sum_{d=0, kd|p-1}^{kf}\binom{kf}{d}a^{kf-d} \bmod p \equiv -\sum_{i=0}^{k/2}\binom{kf}{2fi}a^{f(k-2i)} \bmod p,$$

where the index $i$ runs over the solutions of $kd = (p-1)i = 2kfi$, and $0 \leq d \leq kf$. ∎

**Theorem 14.** Let $k \geq 1$ and let $p = 2kf + 1$ be a prime, then

$$\phi_k(a) \equiv \sum_{x=0}^{p-1}\left(\frac{x}{p}\right)\left(\frac{x^k+a}{p}\right) \equiv -\sum_{i=0}^{(k+1)/2}\binom{kf}{(2i-1)f}a^{f(k-2i+1)} \bmod p. \tag{21}$$

Proof: Routine calculations show that

$$\phi_k(a) = \sum_{x=0}^{p-1}\left(\frac{x}{p}\right)\left(\frac{x^k+a}{p}\right) \equiv \sum_{x=0}^{p-1}(x^{k+1}+ax)^{(p-1)/2} \bmod p \equiv \sum_{d=0}^{kf}\binom{kf}{d}a^{kf-d}\sum_{x=0}^{p-1}x^{k(f+d)} \bmod p. \tag{22}$$

The inner sum vanishes for $k(f+d) \not\equiv 0 \bmod p - 1$, see Lemma 2, but equals $-1$ for $k(f+d) \equiv 0 \bmod p - 1$. Hence

$$\phi_k(a) \equiv -\sum_{d=0,\, p-1|k(f+d)}^{kf}\binom{kf}{d}a^{kf-d} \bmod p \equiv -\sum_{i=0}^{(k+1)/2}\binom{kf}{(2i-1)f}a^{f(k-2i+1)} \bmod p, \tag{23}$$

where the index $i$ runs over the solutions of $k(f+d) = (p-1)i = 2kfi$, and $0 \leq d \leq kf$. ∎

**Note:** The binomial coefficient vanishes in the two cases: $\binom{n}{-k} = 0$ for $k > 0$, and $\binom{n}{k} = 0$ for $k > n \geq 1$.

**Lemma 15.** Let $C : y^2 = x^k + a$, $k > 4$, Then
(i) The group of automorphism $\text{Aut}(C) = \mu_k \times \mu_2$, and defined by $(x, y) \rightarrow (\omega x, \pm y)$, where $\omega \in \mu_k$ is the group of $k$th root of unity.
(ii) The algebraic curve $C$ has complex multiplication by $G = \mu_k \times \mu_2$, and defined by $(x, y) \rightarrow (\omega x, \pm y)$, where $\omega \in \mu_k$ is the group of $k$th root of unity.

## 3. Closed Form Evaluations

The method for evaluating the exponential sums $S(g_n(x))$ in closed form uses the values of simpler exponential quadratic exponential sum $S(ax^2 + bx + c)$ and cubic exponential sum $S(f_n(x))$ of lower degrees to derive the values of the exponential sums $S(g_n(x))$ of higher degree. A different but related method was employed in [WS] to evaluate several classes of quartic exponential sums, including $S(g_n(x))$, $n = 1, 2$, and 7.

### 3.1 The Proof of Theorem 1.

**Theorem 1.** (1) If $n = 1, 2$, or 7, then

$$S(g_n) = \sum_{x=0}^{p-1} \left( \frac{g_n(x)}{p} \right) = \begin{cases} A & \text{if } p \text{ is inert,} \\ A + \left( \frac{a}{p} \right) u & \text{if } 4p = u^2 + nv^2, \end{cases} \tag{24}$$

where the integer $u$ is uniquely determined by quadratic symbol relation $(u \mid n) = (2 \mid p)$, and the quadratic exponential sum

$$A = \sum_{x=0}^{p-1} \left( \frac{f_n(x)/x}{p} \right). \tag{25}$$

(2) If $n = 3, 11, 19, 43, 67$, or 163, then

$$S(g_n) = \sum_{x=0}^{p-1} \left( \frac{g_n(x)}{p} \right) = \begin{cases} B & \text{if } p \text{ is inert,} \\ B + \left( \frac{a}{p} \right) u & \text{if } 4p = u^2 + nv^2, \end{cases} \tag{26}$$

where the integer $u$ is uniquely determined by $(u \mid n) = (2 \mid p)$, $g_n(x) = (x + a_1)(x + a_2)(x + a_3)(x + a_4) \in \mathbf{F}_p[x]$, and

$$B = -1 - \left( \frac{\alpha}{p} \right) (H(\beta) \bmod p) \quad \text{with} \quad \alpha \equiv (a_1 - a_4)(a_2 - a_3), \ \beta \equiv \frac{(a_1 - a_3)(a_2 - a_4)}{(a_1 - a_4)(a_2 - a_3)} \bmod p, \tag{27}$$

and $H(x)$ is the Hasse invariant.

Proof: To verify the first case, it is sufficient to consider $f_1(x) = x^3 + ax$, and $g_1(x) = f_1(x^2) = x^4 + a$, the singular part is removed since the quadratic symbol annihilates it. Apply the transformation formula, Lemma 5, to $f_1(x)$ to obtain

$$\sum_{x=0}^{p-1}\left(\frac{x^6+ax^2}{p}\right) = \sum_{x=0}^{p-1}\left(\frac{x(x^3+ax)}{p}\right) + \sum_{x=0}^{p-1}\left(\frac{x^3+ax}{p}\right). \tag{28}$$

Now, substituting the known values of the exponential sums $A = S(f_1(x)/x)$, see Lemma 3, and $S(f_1(x))$, see Theorem 4, and simplifying return

$$\sum_{x=0}^{p-1}\left(\frac{x^4+a}{p}\right) = \sum_{x=0}^{p-1}\left(\frac{x^2+a}{p}\right) + \sum_{x=0}^{p-1}\left(\frac{x^3+ax}{p}\right) = \begin{cases} A & \text{if } p \text{ is inert,} \\ A + \left(\dfrac{a}{p}\right)u & \text{if } 4p = u^2 + nv^2. \end{cases} \tag{29}$$

To verify the second case, it is sufficient to consider the polynomials $f_{11}(x) = x^3 - 2^5 \cdot 3 \cdot 11 a^2 x + 2^4 \cdot 7 \cdot 11^2 a^3$, and $g_{11}(x) = f_{11}(x^2)$. Apply the transformation formula, Lemma 5, to $f_{11}(x)$ to obtain

$$\sum_{x=0}^{p-1}\left(\frac{f_{11}(x^2)}{p}\right) = \sum_{x=0}^{p-1}\left(\frac{xf_{11}(x)}{p}\right) + \sum_{x=0}^{p-1}\left(\frac{f_{11}(x)}{p}\right). \tag{30}$$

To apply the reduction formula, Lemma 6, factor the polynomial as

$$xf_{11}(x) = x(x^3 - 2^5 \cdot 3 \cdot 11 a^2 x + 2^4 \cdot 7 \cdot 11^2 a^3) = (x+a_1)(x+a_2)(x+a_3)(x+a_4) \in \mathbf{F}_p[x], \tag{31}$$

where $a_1 = 0$, and put $\alpha \equiv (a_1 - a_4)(a_3 - a_2)$, $\beta \equiv ((a_1 - a_4)(a_2 - a_3))^{-1}(a_1 - a_3)(a_2 - a_4) \bmod p$. Next, substitute the known values of the cubic exponential sum $S(f_{11}(x))$, see Theorem 4, and the quartic exponential sum $S(xf_{11}(x))$, Lemmas 6 and 9, to arrive at

$$\sum_{x=0}^{p-1}\left(\frac{f_{11}(x^2)}{p}\right) = \sum_{x=0}^{p-1}\left(\frac{xf_{11}(x)}{p}\right) + \sum_{x=0}^{p-1}\left(\frac{f_{11}(x)}{p}\right) = \begin{cases} B & \text{if } p \text{ is inert,} \\ B + \left(\dfrac{a}{p}\right)u & \text{if } 4p = u^2 + nv^2, \end{cases} \tag{32}$$

where

$$B = \sum_{x=0}^{p-1}\left(\frac{xf_{11}(x)}{p}\right) = -1 - \left(\frac{\alpha}{p}\right)(H(\beta) \bmod p), \tag{33}$$

and the polynomial $H(x)$ is the Hasse invariant, see Lemma 9. ∎

**Note**: The last formula works over the splitting field of a polynomial $f(x) = (x+a_1)(x+a_2)(x+a_3)(x+a_4)$. If it does not split over the prime field $\mathbf{F}_p$, then splitting field of $f(x)$ is either $\mathbf{F}_{p^2}$, $\mathbf{F}_{p^3}$ or $\mathbf{F}_{p^4}$.

The new polynomials $g_n(x)$ are connected with algebraic curves of genus $g = 1, 2$. These curves are as follows.

(1) The elliptic curve $C_n : y^2 = g_n(x)$, with the parameters $n = 1, 2, 7$, and $g = 1$.

(2) The hyperelliptic curve $C_n : y^2 = g_n(x)$, with the parameters $n = 3, 11, 19, 43, 67, 163$, and $g = 2$.

***Corollary* 16.** The number of points on the hyperelliptic curve $C_n : y^2 = g_n(x)$ is $\# C_n(\mathbf{F}_p) = p + 1 + S(g_n(x))$.

***Corollary* 17.** Let $g_n(x) = f_n(x^2)$, then

$$\left| \sum_{x=0}^{p-1} \left( \frac{g_n(x)}{p} \right) \right| \leq 2\sqrt{p}. \tag{34}$$

Proof: Since $H(\beta) \bmod p \equiv a_p < 2g\sqrt{p}$, where $g > 0$ is the genus of the curve, the absolute value of $H(d) \bmod p$ is the same as the absolute value of the trace of Frobenius $a_p = -S(F_d)$ at $p$. ∎

This improves the Weil estimate $|S(f)| \leq (n-1)\sqrt{p}$, $4 \leq n \leq 6$, for the polynomials $f(x) = g_n(x)$, and their twisted polynomials $g(x) = sg_n(x)$, $s \in \mathbf{F}_p$.

**3.2. More Closed Form Evaluations.**
Other widely used families of nonCM elliptic curves of genus $g = 1$ are the followings:

(10) $F_\beta(x) = x(x - 1)(x - \beta)$, or $F_\beta(x) = (x^2 - 1)(x^2 - \beta)$, where $0, 1 \neq \beta \in \mathbb{Z}$,      Legendre form.
(11) $G_\beta(x) = (k^2x^2 - 1)(x^2 - \beta)$, $0, \pm 1 \neq \beta \in \mathbb{Z}$, and $k \neq 0$,      Newton form.
(12) $H_d(x) = (x^2 - c^2)(c^2dx^2 - 1)$, $c, d \in \mathbb{Z}$, and $cd(1 - c^4d) \neq 0$,      Edward form.

The Edward form $H_d(x) = (x^2 - c^2)(c^2dx^2 - 1)$ is derived from the elliptic curve $E: x^2 + y^2 = c^2(1 + dx^2y^2)$ or its quadratic twist $E': sx^2 + y^2 = 1 + dx^2y^2$, see [ES], and [BL]. Basically, the Newton form $E: y^2 = x^4 - ax^2 + 1$, and Jacobi elliptic curve $E: y^2 = x^4 - ax^2 + 1$ are special cases of the Edward form.

Since the Edward form $H_\beta(x) = (x^2 - c^2)(c^2dx^2 - 1) = (x + a_1)(x + a_2)(x + a_3)(x + a_4) \in \mathbf{F}_p[x]$ or over a quadratic extension $\mathbf{F}_q$ of $\mathbf{F}_p$, where $a_1 = c$, $a_2 = -c$, $a_3 = -1/c\sqrt{d}$, $a_4 = 1/c\sqrt{d}$, an elliptic curve $E: y^2 = x^3 + ax + b$, with 2-torsion points is isogenous to an Edward form.

Similarly, an arbitrary elliptic curve $E: y^2 = f(x)$ is isogenous to the Legendre form $E': y^2 = F_\beta(x)$, or one of the other forms over some finite extension $\mathbf{F}_q$ of $\mathbf{F}_p$, it appears that the associated exponential sum $S(f)$ can be evaluated in closed form in term of simpler exponential sums such as $S(F_\beta)$ and quadratic exponential sums $S(ax^2 + bx + c)$. A new paper, see [OG], develops the rational maps between the Legendre forms and the Edward forms of elliptic curves.

***Theorem* 18.** Let $F(x) = F_\beta(x)$, $G_\beta(x)$, or $H_d(x)$, then

$$S(F) = \sum_{x=0}^{p-1} \left( \frac{F(x)}{p} \right) = -1 - \left( \frac{\alpha}{p} \right) (H(\beta) \bmod p) \tag{35}$$

where $-a_1, -a_2, -a_3, -a_4$ are the roots of the polynomial $F(x)$, and $\alpha \equiv (a_1 - a_3)(a_3 - a_2)$, and $\beta \equiv ((a_1 - a_4)(a_2 - a_3))^{-1} a_3(a_2 - a_4) \bmod p$.
Proof: Factor the polynomial as $F(x) = (x + a_1)(x + a_2)(x + a_3)(x + a_4) \in \mathbf{F}_q[x]$, and put

$$\alpha \equiv (a_1 - a_4)(a_2 - a_3), \quad \text{and} \quad \beta \equiv \frac{(a_1 - a_3)(a_2 - a_4)}{(a_1 - a_4)(a_2 - a_3)} \bmod p. \tag{36}$$

By Lemmas 6 and 9,

$$S(F) = \sum_{x=0}^{p-1} \left(\frac{F(x)}{p}\right) = -1 - \left(\frac{\alpha}{p}\right) S(F_\beta). \tag{37}$$

Next, substitute the known value of the exponential sum $S(F_\beta(x)) \equiv H(\beta) \bmod p$, Lemma 9, and simplify. The evaluations are correct up to a permutation of the roots $a_{\pi(1)}, a_{\pi(2)}, a_{\pi(3)}, a_{\pi(4)}$, where $\pi \in \text{Sym}(4)$ is a permutation of 4 elements, of the roots of the polynomial $F(x)$. ∎

The exponential sum $S(G_\beta(x))$ for the special case $G_\beta(x) = (k^2 x^2 - 1)(x^2 - \beta)$, $0, \pm 1 \neq \beta \in \mathbb{Z}$, for $k = 1$, can be computed in a slightly different way as

$$\sum_{x=0}^{p-1}\left(\frac{(x^2-1)(x^2-\beta)}{p}\right) = \sum_{x=0}^{p-1}\left(\frac{x(x-1)(x-\beta)}{p}\right) + \sum_{x=0}^{p-1}\left(\frac{(x-1)(x-\beta)}{p}\right). \tag{38}$$

Completing the calculations yield

$$\sum_{x=0}^{p-1}\left(\frac{(x^2-1)(x^2-\beta)}{p}\right) = A_\beta + (H(\beta) \bmod p) \tag{39}$$

where

$$A_\beta = \sum_{x=0}^{p-1}\left(\frac{(x-1)(x-\beta)}{p}\right) = \begin{cases} -1 & \text{if } (1-\beta)^2 \neq 0, \\ p-1 & \text{if } (1-\beta)^2 = 0, \end{cases} \tag{40}$$

see Lemma 3.

**Theorem 19.** Let $F(x) = x^k - a$, $a \in \mathbb{F}_q$. If $\gcd(p-1, k) \neq 1$, then

$$\sum_{x=0}^{p-1}\left(\frac{x^{2k}-a}{p}\right) = -\sum_{i=0}^{(k+1)/2}\binom{kf}{(2i-1)f} a^{f(k-2i+1)} + \binom{kf}{2fi} a^{f(k-2i)} \bmod p. \tag{41}$$

Proof: Using the identity

$$\sum_{x=0}^{p-1}\left(\frac{g_n(x^2)}{p}\right) = \sum_{x=0}^{p-1}\left(\frac{xg_n(x)}{p}\right) + \sum_{x=0}^{p-1}\left(\frac{g_n(x)}{p}\right). \tag{42}$$

In full details, this is precisely

$$\sum_{x=0}^{p-1}\left(\frac{x^{2k}+a}{p}\right) = \sum_{x=0}^{p-1}\left(\frac{x(x^2+a)}{p}\right) + \sum_{x=0}^{p-1}\left(\frac{x^2+a}{p}\right) = \phi_k(a) + \psi_k(a), \tag{43}$$

where the last two exponential sums are, see Theorem 13, and 14,

$$\phi_k(a) \equiv \sum_{x=0}^{p-1} \left(\frac{x}{p}\right)\left(\frac{x^k+a}{p}\right) \equiv -\sum_{i=0}^{(k+1)/2} \binom{kf}{(2i-1)f} a^{f(k-2i+1)} \bmod p, \tag{44}$$

and

$$\psi_k(a) = \sum_{x=0}^{p-1} \left(\frac{x^k+a}{p}\right) \equiv -\sum_{i=0}^{k/2} \binom{kf}{2fi} a^{f(k-2i)} \bmod p \text{ if } \gcd(p-1, k) \neq 1. \tag{45}$$